\documentclass[a4paper,10pt]{article}

\usepackage[ascii]{inputenc}
\usepackage{amsmath,amsfonts,amssymb,amsthm}

\newtheoremstyle{theorem}{10pt}{5pt}{\slshape}{}{\bfseries}{.}{.5em}{}

\theoremstyle{theorem}
\newtheorem{theorem}{Theorem}
\newtheorem{lemma}[theorem]{Lemma}
\newtheorem{corollary}[theorem]{Corollary}

\usepackage{titlesec}
\titleformat{\section}
{\normalfont\bfseries}{\thesection}{1em}{\large}
\titlespacing*{\section}{0pt}{3.5ex plus 1ex minus .2ex}{2.3ex plus .2ex}

\titleformat{\subsection}
{\normalfont\bfseries}{\thesubsection}{1em}{\normalsize}
\titlespacing*{\section}{0pt}{3.5ex plus 1ex minus .2ex}{2.3ex plus .2ex}

\allowdisplaybreaks[3]

\begin{document}

\title{Rellich Type Theorems for Unbounded Domains}
\author{Esa V\!. Vesalainen}
\date{}
\maketitle

\begin{abstract}
We give several generalizations of Rellich's classical uniqueness theorem to unbounded domains. We give a natural half-space generalization for super-exponentially decaying inhomogeneities using real variable techniques. We also prove under super-exponential decay a discrete generalization where the inhomogeneity only needs to vanish in a suitable cone.

The more traditional complex variable techniques are used to prove the half-space result again, but with less exponential decay, and a variant with polynomial decay, but with supports exponentially thin at infinity. As an application, we prove the discreteness of non-scattering energies for non-compactly supported potentials with suitable asymptotic behaviours and supports.
\end{abstract}

\section{Introduction}

\subsection{Scattering theory}

Our objects of study arise from scattering theory. More precisely, time independent scattering theory for short-range potentials, which models e.g. two-body quantum scattering, acoustic scattering, and some classical electromagnetic scattering situations (for a general reference, see e.g. \cite{Colton--Kress}).
Here one is concerned with the situation where, at a fixed energy or wavenumber $\lambda\in\mathbb R_+$, an incoming wave $w$, which is a solution to the free equation \[\left(-\Delta-\lambda\right)w=0,\] is scattered by some perturbation of the flat homogeneous background. Here this perturbation will be modeled by a real-valued function $V$ in $\mathbb R^n$ having enough decay at infinity. The total wave $v$, which models the ``actual'' wave, then solves the perturbed equation \[\left(-\Delta+V-\lambda\right)v=0.\] For acoustic and electromagnetic scattering, one writes $\lambda V$ instead of $V$.
Of course, the two waves $v$ and $w$ must be linked together and the connection is given by the Sommerfeld radiation condition. The upshot will be that the difference $u$ of $v$ and $w$, the so-called scattered wave, will have an asymptotic expansion of the shape
\[u(x)=v(x)-w(x)=A\!\left(\frac x{\left|x\right|}\right)\frac{e^{i\sqrt\lambda\left|x\right|}}{\left|x\right|^{(n-1)/2}}+\text{error},\]
where $A$ depends on $\lambda$ and $w$, and where the error term decays more rapidly than the main term. The point here is that in the main term the dependences on the radial and angular variables are neatly separated, and in practical applications one usually measures the scattering amplitude or far-field pattern $A$, or its absolute value $\left|A\right|$.

\subsection{Non-scattering energies}

%An example of a stereotypical question in scattering theory would be: what can be said about the potential $V$ if $A$ is known for many $w$?
%Two celebrated methods, the linear sampling method \cite{} and the factorization method \cite{}, can be used to derive information about the support of $V$, when the support is compact. However, the energy $\lambda$ needs to be such that $A\not\equiv0$ for every $w\not\equiv0$.

It is a natural question whether we can have $A\equiv0$ for some $w\not\equiv0$? This would mean that the main term of the scattered wave vanishes at infinity, meaning that the perturbation, for the special incident wave in question, is not seen far away. Values of $\lambda\in\mathbb R_+$ for which such an incident wave $w$ exists, are called non-scattering energies (or appropriately, wavenumbers) of $V$. In order to avoid the discussion of function spaces here, the precise definition is given in Section \ref{main-result-section} below.

Results on the existence of non-scattering energies are scarce. Essentially only two general results are known: For compactly supported radial potentials the set of non-scattering energies is an infinite discrete set accumulating at infinity \cite{Colton--Monk}, and for compactly supported potentials with suitable corners, Bl{\aa}sten, P\"aiv\"arinta and Sylvester recently proved that the set of non-scattering energies is empty~\cite{Blasten--Paivarinta--Sylvester}.

We would like to mention the related topic of transparent potentials: there one considers (at a fixed energy) potentials for which $A$ vanishes for all $w$. The knowledge of transparent potentials is more extensive. In particular, several constructions of such radial potentials have been given, see e.g. the works of Regge \cite{Regge}, Newton \cite{Newton}, Sabatier \cite{Sabatier}, Grinevich and Manakov \cite{Grinevich--Manakov}, and Grinevich and Novikov \cite{Grinevich--Novikov}.

\subsection{Rellich type theorems}

In practice, discreteness of the set of non-scattering energies tends to be a more attainable goal. The first key step towards that goal (for compactly supported~$V$) is supplied by Rellich's classical uniqueness theorem which is the following:
\begin{theorem}\label{rellich-lemma}
Let $u\in L^2_{\mathrm{loc}}(\mathbb R^n)$ solve the equation $\left(-\Delta-\lambda\right)u=f$, where $\lambda\in\mathbb R_+$ and $f\in L^2(\mathbb R^n)$ is compactly supported, and assume that
\[\frac1R\int\limits_{B(0,R)}\left|u\right|^2\longrightarrow0,\]
as $R\longrightarrow\infty$.
Then $u$ also is compactly supported.
\end{theorem}

This was first proved (though with a bit different decay condition) independently by Rellich \cite{Rellich} and Vekua \cite{Vekua} in 1943.
Of the succeeding work, which includes generalizations of this result to more general constant coefficient differential operators, we would like to mention the work of Tr\`eves \cite{Treves}, Littman \cite{Littman1, Littman2, Littman3}, Murata \cite{Murata} and H\"ormander \cite{Hormander}. Section 8 of \cite{Hitrik--Krupchyk--Ola--Paivarinta} is also interesting.

We also mention that a theorem analogous to Theorem \ref{rellich-lemma} also exists for the discrete Laplacian (also to be defined more precisely in Section \ref{main-result-section}):
\begin{theorem}\label{discrete-rellich-lemma}
Let $u\colon\mathbb Z^n\longrightarrow\mathbb C$ be a solution to $\left(-\Delta_{\mathrm{disc}}-\lambda\right)u=f$, where
\[\frac1R\sum_{\xi\in\mathbb Z^n,\left|\xi\right|\leqslant R}\left|u(\xi)\right|^2\longrightarrow0,\]
as $R\longrightarrow\infty$, and $f\in\ell^2(\mathbb Z^n)$ is non-zero only at finitely many points of $\mathbb Z^n$, and $\lambda\in\left]0,n\right[$.
Then $u$ also is non-zero only at finitely many points of $\mathbb Z^n$.
\end{theorem}
\noindent This theorem was proved recently by Isozaki and Morioka \cite{Isozaki--Morioka}. A less general version of the result was implicit in the work of Shaban and Vainberg \cite{Shaban--Vainberg}.

\subsection{Transmission eigenvalues}

Assume that $V$ is compactly supported. The equations for $v$ and $w$ imply that the scattered wave $u$ solves the equation
\[\left(-\Delta-\lambda\right)u=-Vv.\]
If furthermore $A\equiv0$, then $u$ will satisfy the decay condition in Theorem \ref{rellich-lemma}, and so $u=v-w$ will vanish outside a compact set. If the support of $V$ is essentially some suitable open domain $\Omega$, the unique continuation principle for the free Helmholtz equation allows us to conclude that actually
\[\left\{\!\!\begin{array}{ll}
\left(-\Delta+V-\lambda\right)v=0&\text{in $\Omega$,}\\
\left(-\Delta-\lambda\right)w=0&\text{in $\Omega$,}\\
v-w\in H_0^2(\Omega).\end{array}\right.\]
This system, called the interior transmission problem, is a non-self-adjoint eigenvalue problem for $\lambda$, and the values of $\lambda$, for which this system has non-trivial $L^2$-solutions, are called (interior) transmission eigenvalues.

The non-scattering energies and transmission eigenvalues first appeared in the papers of Colton and Monk \cite{Colton--Monk} and Kirsch \cite{Kirsch}. In \cite{Colton--Kirsch--Paivarinta} Colton, Kirsch and P\"aiv\"arinta proved the discreteness of transmission eigenvalues (and non-scattering energies) for potentials that may even be mildly degenerate. The early papers on the topic also considered, among other things, radial potentials; for more on this, we refer to the article of Colton, P\"aiv\"arinta and Sylvester~\cite{Colton--Paivarinta--Sylvester}.

In recent years, there has been a surge of interest in the topic starting with the general existence results of P\"aiv\"arinta and Sylvester \cite{Paivarinta--Sylvester}, who established existence of transmission eigenvalues for a large class of potentials, and Cakoni, Gintides and Haddar \cite{Cakoni--Gintides--Haddar}, who established for acoustic scattering, that actually the set of transmission eigenvalues must be infinite.

For potentials more general than the radial ones, a very common approach to proving discreteness and other properties has been via quadratic forms: the scattered wave solves the fourth-order equation
\[\left(-\Delta+V-\lambda\right)\frac1V\left(-\Delta-\lambda\right)u=0,\]
and this can be handled nicely with quadratic forms (or with variational formulations) and analytic perturbation theory.

Recently, other approaches, not involving the fourth-order equation, to proving discreteness and many other results have been introduced by Sylvester \cite{Sylvester}, Robbiano \cite{Robbiano}, and Lakshtanov and Vainberg \cite{Lakshtanov--Vainberg}.

For more information and a wealth of references on transmission eigenvalues, we recommend the survey of Cakoni and Haddar \cite{Cakoni--Haddar} and their editorial \cite{Cakoni--Haddar2}.

\subsection{What we do and why?}

Most of the work on non-scattering energies and transmission eigenvalues deals with compactly supported potentials $V$. However, the basic short-range scattering theory only requires $V$ to have enough decay at infinity, essentially something like $V(x)\ll\left|x\right|^{-1-\varepsilon}$. Thus, it makes perfect sense to study non-scattering energies for non-compactly supported potentials.

In \cite{Vesalainen}, we studied an analogue of the transmission eigenvalue problem for certain unbounded domains establishing discreteness under the assumptions that $V(x)\asymp\left|x\right|^{-\alpha}$ for some $\alpha\in\left]3,\infty\right[$, and that the underlying domain $\Omega$ has the property that the embedding $H^2_0(\Omega)\longrightarrow L^2(\Omega)$ is compact.

Even though the formulation of the transmission problem was kept compatible with scattering theory, the lack of a Rellich type theorem for unbounded domains did not allow any conclusions about the discreteness of non-scattering energies for the corresponding potentials. Here we will present several such Rellich type theorems, two of which will give discreteness of non-scattering energies.

First, we prove a fairly general Rellich type theorem where instead of a compactly supported inhomogeneity $f$, we consider $f$ that is superexponentially decaying and vanishes in a half-space. The conclusion will then be that the solution $u$ also vanishes in the same half-space. This is a fairly satisfying generalization. Also, it contains the classical Rellich lemma as a simple corollary. We give a new kind of proof for this result based on real variable techniques, first deriving a Carleman estimate weighted exponentially in one direction from an  estimate of Sylvester and Uhlmann \cite{Sylvester--Uhlmann} and then arguing immediately from it.

However, we have so far been unable to apply the quadratic form approach to superexponentially decaying potentials, and so we would like to have Rellich type theorems which allow less decay. We shall give two results of this kind: the first is for exponentially decaying inhomogeneities, the second is for polynomially decaying potentials but for domains that are not only contained in a half-space but also grow exponentially thin at infinity. These results are proved with a more traditional complex variables argument \cite{Treves, Littman1, Littman2, Littman3, Hormander}.

The Rellich type theorem for polynomially decaying inhomogeneities can be immediately combined with the results of \cite{Vesalainen} to give discreteness of non-scattering energies for a class of polynomially decaying potentials. Obtaining a discreteness result for a class of exponentially decaying potentials will require some minor adjustments to the arguments of \cite{Vesalainen} which are presented in Section \ref{corollary-section}.
Our manner of using quadratic forms to establish discreteness is a close relative of the application of quadratic forms to degenerate and singular potentials in the works of Colton, Kirsch and P\"aiv\"arinta \cite{Colton--Kirsch--Paivarinta}, Serov and Sylvester\cite{Serov--Sylvester}, Serov \cite{Serov}, and Hickmann \cite{Hickmann}.

Finally, as an interesting aside, and to provide a point of comparison, we present a generalization of the discrete Rellich type theorem of Isozaki and Morioka. It turns out that for superexponentially decaying potentials, one gets a much stronger result than in the continuous case: we not only can consider vanishing in half-spaces but vanishing in suitable cones. The proof depends heavily on the arguments in \cite{Isozaki--Morioka} which are first used to show that the solution must be superexponentially decaying. After this, the Rellich type conclusion follows from a repeated application of the definition of the discrete Laplacian.

\subsection{On notation}

We shall use the standard asymptotic notation. If $f$ and $g$ are complex functions on some set $A$, then $f\ll g$ means that $\left|f(x)\right|\leqslant C\left|g(x)\right|$ for all $x\in A$ for some positive real constant $C$, referred to as the implicit constant. The notation $f\asymp g$ means that both $f\ll g$ and $g\ll f$. When $C$ may depend on some parameters $\alpha$, $\beta$, \dots, we write $f\ll_{\alpha,\beta,\ldots}g$, except that all the implicit constants are allowed to depend on the dimension $n\in\left\{2,3,\ldots\right\}$ of the ambient Euclidean space $\mathbb R^n$.

For a vector $x\in\mathbb R^n$, we define $\left\langle x\right\rangle=\sqrt{1+\left|x\right|^2}$ and $x'=\left\langle x_1,\ldots,x_{n-1}\right\rangle$. The letters $e_1$, $e_2$, \ldots, $e_n$ denote the standard basis of $\mathbb R^n$:
\begin{align*}
e_1&=\left\langle1,0,0,\ldots,0,0\right\rangle,\\
e_2&=\left\langle0,1,0,\ldots,0,0\right\rangle,\\
&.........................\\
e_n&=\left\langle0,0,0,\ldots,0,1\right\rangle.
\end{align*}
For a complex vector $z\in\mathbb C^n$, we denote the real and imaginary parts of $z$ by
\[\Re z=\left\langle\Re z_1,\ldots,\Re z_n\right\rangle\quad\text{and}\quad
\Im z=\left\langle\Im z_1,\ldots,\Im z_n\right\rangle.\]
For $R\in\mathbb R_+$, we write $B(0,R)$ for the open ball of vectors $\xi\in\mathbb R^n$ with $\left|\xi\right|<R$. If we want to emphasize the dimension of the ambient Euclidean space, we write $B^n(0,R)$.

We shall use the shorthand $\nabla\otimes\nabla$ to simplify expressions of Sobolev norms in the obvious way. For example, the usual $H^2$-Sobolev norm would be given by the expression
\[\sqrt{\bigl\|u\bigr\|^2+\bigl\|\nabla u\bigr\|^2+\bigl\|\nabla\otimes\nabla u\bigr\|^2}.\]

We use $\mathbb T^n$ to denote the $n$-dimensional torus $\mathbb R^n/\mathbb Z^n$. The corresponding complex torus $\mathbb T^n_{\mathbb C}$ means $\mathbb C^n/\mathbb Z^n$. Instead of analytic functions in $\mathbb T^n_{\mathbb C}$ one can simply think of entire functions in $\mathbb C^n$ which are $1$-periodic with respect to each complex variable.

For a function $f\in L^2(\mathbb R^n)$, $\widehat f$ denotes the usual Fourier transform normalized as follows: for a Schwartz test function $f$ and $\xi\in\mathbb R^n$,
\[\widehat f(\xi)=\int\limits_{\mathbb R^n}f(x)\,e(-x\cdot\xi)\,\mathrm dx,\]
where $e(\cdot)$ stands for $e^{2\pi i\cdot}$. We will also denote by $F'f$ the Fourier transform of $f$ in the first $n-1$ variables, and the corresponding inverse transform by $F'^{-1}$. Similarly, we will denote by $F_n$ the Fourier transform in the $n$th variable, and by $F_n^{-1}$ the corresponding inverse transform.

In the same vein, given a vector $f\in\ell^2(\mathbb Z^n)$, $\check f$ denotes the Fourier series connected to $f$: for $x\in\mathbb T^n$,
\[\check f(x)=\sum_{\xi\in\mathbb Z^n}f(x)\,e(x\cdot\xi),\]
and the convergence is in the $L^2$-sense or pointwise, whichever is more appropriate.

Finally, the usual $L^2$-inner product is denoted by $\left\langle\cdot\middle|\cdot\right\rangle$: for square-integrable complex functions $f$ and $g$ in some domain $\Omega$, we define
\[\left\langle f\middle|g\right\rangle=\int\limits_\Omega\overline f\,g.\]

\section{The results}\label{main-result-section}

\subsection{Rellich type theorems for unbounded domains}

Before stating the main results, we would like to define the function spaces from which solutions $v$, $w$ and $u$ are actually taken. The solutions $v$ and $w$ should be taken from the Agmon--H\"ormander space $B^\ast$, which consists of those $L_{\mathrm{loc}}^2(\mathbb R^n)$-functions $u$ for which
\[\left\|u\right\|_{B^\ast}^2=\sup_{R>1}\frac1R\int\limits_{B(0,R)}\left|u\right|^2<\infty.\]
Actually, it will then turn out that also the first- and second-order partial derivatives of $v$ and $w$ will also belong to $B^\ast$, a fact which we will denote by $v,w\in B^\ast_2$.

When $A\equiv0$, the scattered wave $u$ will satisfy
\[\frac1R\int\limits_{B(0,R)}\left|u\right|^2\longrightarrow0\]
as $R\longrightarrow\infty$. The space of such $L_{\mathrm{loc}}^2(\mathbb R^n)$-functions is denoted by $\mathring B^\ast$. From Theorem 14.3.6 of \cite{Hormander-II}, it will follow that the first- and second-order partial derivatives of $u$ will also belong to $\mathring B^\ast$, and we will write $u\in\mathring B^\ast_2$. In view of this, the growth condition for $u$ in Theorem \ref{rellich-lemma} could be replaced by $u\in\mathring B^\ast_2$, and this is what we shall do in the theorems below.

The Agmon--H\"ormander spaces were introduced in the works of Agmon and H\"ormander \cite{Agmon--Hormander}, and independently by Murata \cite{Murata-1}, in the study of asymptotics of solutions to constant coefficient partial differential equations. Despite their appearance, these spaces turn out to be fairly natural. In particular, the mapping from $g\in L^2(S^{n-1})$ to the Herglotz wave
\[w(x)=\int_{S^{n-1}}g(\vartheta)\,e^{i\sqrt\lambda\vartheta\cdot x}\,\mathrm d\vartheta\]
will give a bijection from $L^2(S^{n-1})$ into the $B^\ast$-solutions of $\left(-\Delta-\lambda\right)w=0$, and the respective norms of $g$ and the Herglotz wave $w$ will be comparable.

Our first generalization of Theorem \ref{rellich-lemma} considers vanishing in a half-space instead of the exterior of a ball.
\begin{theorem}\label{superexponential-main-theorem}
Let $u\in\mathring B^\ast_2$ solve
\[\left(-\Delta-\lambda\right)u=f,\]
where $\lambda\in\mathbb R_+$ and $f\in e^{-\gamma\left\langle\cdot\right\rangle}L^2(\mathbb R^n)$ for all $\gamma\in\mathbb R_+$, and suppose that $f$ vanishes in the lower half-space $\mathbb R^{n-1}\times\mathbb R_-$. Then also $u$ vanishes in $\mathbb R^{n-1}\times\mathbb R_-$.
\end{theorem}
With complex variable techniques, we may allow less exponential decay:
\begin{theorem}\label{exponential-main-theorem}
Let $u\in\mathring B_2^\ast$ solve
\[\left(-\Delta-\lambda\right)u=f,\]
where $\lambda\in\mathbb R_+$ and $f\in e^{-\gamma_0\left\langle\cdot\right\rangle}L^2(\mathbb R^n)$ for some $\gamma_0\in\mathbb R_+$, and suppose that $f$ vanishes in $\mathbb R^{n-1}\times\mathbb R_-$.
Then also $u$ vanishes in $\mathbb R^{n-1}\times\mathbb R_-$.
\end{theorem}

Finally, for polynomially decaying $f$ we assume that the support of $f$ is exponentially thin at infinity:
\begin{theorem}\label{polynomial-main-theorem}
Let $u\in\mathring B_2^\ast$ solve
\[\left(-\Delta-\lambda\right)u=f,\]
where $\lambda\in\mathbb R_+$ and the inhomogeneity $f$ both vanishes in $\mathbb R^{n-1}\times\mathbb R_-$ and satisfies
\[\sup_{\substack{\zeta\in\mathbb C^{n-1}\\\left|\Im\zeta\right|<\gamma_0}}\,\,\int\limits_{\mathbb R^n}\left|\left\langle x\right\rangle^2\,f\!\left(x\right)e\!\left(x\cdot\zeta\right)\right|\mathrm dx<\infty\]
for some $\gamma_0\in\mathbb R_+$.
Then $u$ vanishes in the lower half-space $\mathbb R^{n-1}\times\mathbb R_-$.
\end{theorem}

\subsection{Applications to discreteness of non-scattering energies}

The next two theorems will concern discreteness of non-scattering energies in a somewhat special situations. The quadratic form techniques in Section \ref{corollary-section} depend on the potential $V$ being essentially supported in a domain $\Omega$ for which the embedding $H_0^2(\Omega)\longrightarrow L^2(\Omega)$ is compact.

For a discussion of such compact embeddings for unbounded domains, we refer to Chapter 6 of \cite{AdamsEtFournier} or to the original article \cite{Adams}. The conditions get slightly more involved in higher dimensions, but in two- and three-dimensional cases, the embedding $H_0^2(\Omega)\longrightarrow L^2(\Omega)$ is compact if and only if the domain $\Omega$ does not contain an infinite sequence of pairwise disjoint congruent balls; cf. remarks 6.17.3, 6.9 and 6.11 in \cite{AdamsEtFournier}.

Combining Theorem~\ref{polynomial-main-theorem} with the results in \cite{Vesalainen} easily gives the following discreteness result for non-scattering energies:
\begin{theorem}\label{easy-corollary}\setlength{\itemsep}{0pt}
Let $V\in L^\infty\!\left(\mathbb R^n\right)$ take only nonnegative real values, and let $\Omega\subseteq\mathbb R^{n-1}\times\mathbb R_+$ be a non-empty open set for which the Sobolev embedding $H_0^2(\Omega)\longrightarrow L^2(\Omega)$ is compact. Assume the following:
\begin{itemize}\setlength{\itemsep}{0pt}
\item[\rm I.] $V(\cdot)\asymp\left\langle\cdot\right\rangle^{-\alpha}$ in $\Omega$ for some $\alpha\in\left]3,\infty\right[$, and $V$ vanishes in $\mathbb R^n\setminus\Omega$.
\item[\rm II.] The complement of $\Omega$ in $\mathbb R^n$ has a connected interior and is the closure of the interior.
\item[\rm III.] The integrals
\[\int\limits_{\Omega}\left|\left\langle x\right\rangle^{2+1/2+\varepsilon-\alpha}\,e(x\cdot\zeta)\right|^2\,\mathrm dx\]
are uniformly bounded for all $\zeta\in\mathbb C^n$ with $\left|\Im\zeta\right|<\gamma_0$ for some $\gamma_0,\varepsilon\in\mathbb R_+$.
\end{itemize}
Then the set of non-scattering energies for $V$ is a discrete subset of $\left[0,\infty\right[$, and each of them is of finite multiplicity.
\end{theorem}
\noindent
Here multiplicity is defined to be the dimension of the vector space of pairs of solutions $\left\langle v,w\right\rangle$ appearing in the definition of non-scattering energies.
The point of the condition III is that, combined with the Cauchy--Schwarz inequality, and the definition of $B^\ast$, it guarantees that $Vv$ satisfies the $\sup$-condition of Theorem \ref{polynomial-main-theorem} for any $v\in B^\ast$.

Analogously to the Rellich type theorems, if the condition of polynomial decay is replaced by exponential decay, then we can relax the conditions:
\begin{theorem}\label{hard-corollary}\setlength{\itemsep}{0pt}
Let $V\in L^\infty\!\left(\mathbb R^n\right)$ take only nonnegative real values, and let $\Omega\subseteq\mathbb R^{n-1}\times\mathbb R_+$ be a non-empty open set for which the Sobolev embedding $H_0^2(\Omega)\longrightarrow L^2(\Omega)$ is compact. Assume the following:
\begin{itemize}\setlength{\itemsep}{0pt}
\item[\rm I.] $V(\cdot)\asymp e^{-\gamma_0\left\langle\cdot\right\rangle}$ in $\Omega$ for some $\gamma_0\in\mathbb R_+$ with $\gamma_0\ll_n1$, and $V$ vanishes in $\mathbb R^n\setminus\Omega$.
\item[\rm II.] The complement of $\Omega$ in $\mathbb R^n$ has a connected interior and is the closure of the interior.
\end{itemize}
Then the set of non-scattering energies for $V$ is a discrete subset of $\left[0,\infty\right[$ and each of them is of finite multiplicity.
\end{theorem}
\noindent The proof will be similar to the arguments in \cite{Vesalainen}. We will outline the relevant modifications in Section~\ref{corollary-section}.

\subsection{A discrete Rellich type theorem for unbounded domains}

Let $u\in\mathbb Z^n\longrightarrow\mathbb C$ be a function on the square lattice $\mathbb Z^n$. Then we define the discrete Laplacian of $u$ to be the function $-\Delta_{\mathrm{disc}}u\colon\mathbb Z^n\longrightarrow\mathbb C$ given by the formula
\[\bigl(-\Delta_{\mathrm{disc}}u\bigr)(\xi)
=\frac n2u(\xi)-\frac14\sum_{\ell=1}^n\bigl(u(\xi+e_\ell)+u(\xi-e_\ell)\bigr).\]
The spectrum of $-\Delta_{\mathrm{disc}}$ is $\left[0,n\right]$ and absolutely continuous. For more information about the discrete setting, we refer to the presentation of Isozaki and Morioka \cite{Isozaki--Morioka} and the references given there.

Our generalization of Theorem \ref{discrete-rellich-lemma} will concern vanishing, not in a half-space, but in a cone-like domain.
\begin{theorem}\label{discrete-main-theorem}
Let $C$ be the set of those $\xi\in\mathbb Z^n$ for which
\[\left|\xi_1\right|+\left|\xi_2\right|+\ldots+
\left|\xi_{n-1}\right|\leqslant\xi_n.\]
Also, let $u\colon\mathbb Z^n\longrightarrow\mathbb C$ be such that
\[\frac1R\sum_{\left|\xi\right|\leqslant R}\bigl|u(\xi)\bigr|^2\longrightarrow0\]
as $R\longrightarrow\infty$, and let $f\in\ell^2(\mathbb Z^n)$ be such that
\[e^{\gamma\left\langle\cdot\right\rangle}f\in\ell^2(\mathbb Z^n)\] for all $\gamma\in\mathbb R_+$, and assume that $f(\xi)=0$ for all $\xi\in C$.
Finally, let $\lambda\in\left]0,n\right[$, and assume that
\[\left(-\Delta_{\mathrm{disc}}-\lambda\right)u=f\]
in $\mathbb Z^n$. Then also $u(\xi)=0$ for all $\xi\in C$.
\end{theorem}

\section{Proof of Theorem \ref{superexponential-main-theorem} via Carleman estimates}

The following is essentially Lemma 2.5 in \cite[p.~1780]{Paivarinta--Salo--Uhlmann}, which may be reformulated in the following way in view of Theorem 14.3.6 in \cite{Hormander-II}.
\begin{theorem}\label{exponential-rellich}
Let $\gamma_0\in\mathbb R_+$, and let $f$ be a function such that $e^{\gamma\left\langle\cdot\right\rangle}\,f\in L^2(\mathbb R^n)$ for each $\gamma\in\left]0,\gamma_0\right[$. If $u\in\mathring B^\ast_2$ solves the equation
\[\left(-\Delta-\lambda\right)u=f,\]
then
$e^{\gamma\left\langle\cdot\right\rangle}\,u\in H^2(\mathbb R^n)$ for each $\gamma\in\left]0,\gamma_0\right[$.
\end{theorem}

The following result is Proposition 2.1 in \cite{Sylvester--Uhlmann}. Even though the statement there has $\lambda=0$, the same proof also works for $\lambda\geqslant0$.
\begin{theorem}\label{resolvent-estimate}
Let $\lambda\in\mathbb R_+$, $\delta\in\left]-1,0\right[$, and let $\rho\in\mathbb C^n$ with $\rho\cdot\rho=\lambda$ and $\left|\Im\rho\right|\geqslant1$. Then for any $f\in\left\langle\cdot\right\rangle^{-1-\delta}L^2(\mathbb R^n)$ the equation
\[\left(-\Delta-2i\rho\cdot\nabla\right)v=f\]
has a unique solution $v\in\left\langle\cdot\right\rangle^{-\delta}L^2(\mathbb R^n)$ satisfying the estimate
\[\bigl\|\left\langle\cdot\right\rangle^\delta v\bigr\|_{L^2(\mathbb R^n)}
\ll_{\lambda,\delta}\frac1{\left|\rho\right|}\bigl\|\left\langle\cdot\right\rangle^{1+\delta}f\bigr\|_{L^2(\mathbb R^n)}.\]
\end{theorem}
This result can be turned into a Carleman estimate weighted exponentially in one coordinate direction:
\begin{corollary}\label{carleman-estimate}
Let $u\in C_{\mathrm c}^\infty(\mathbb R^n)$, $\delta\in\left]-1,0\right[$, $\lambda\in\mathbb R_+$, and let $\tau\in\mathbb R_+$ with $\tau\gg_{\lambda,\delta}1$. Then
\[\bigl\|e^{-\tau x_n}\left\langle\cdot\right\rangle^\delta u\bigr\|_{L^2(\mathbb R^n)}
\ll_{\lambda,\delta}\frac1\tau\bigl\|e^{-\tau x_n}\left\langle\cdot\right\rangle^{\delta+1}\left(-\Delta-\lambda\right)u\bigr\|_{L^2(\mathbb R^n)}.\]
\end{corollary}
Of course, this also holds for $u$ in the closure of $C_{\mathrm c}^\infty(\mathbb R^n)$ in the norm that is given by the sum of the norms appearing in the estimate.

\paragraph{Proof.} We shall apply Theorem \ref{resolvent-estimate} with
$\rho=\alpha-i\tau e_n$, where $e_n=\left\langle0,0,\ldots,0,1\right\rangle$ and $\alpha\in\mathbb R^n$ is such that $\alpha\cdot e_n=0$ and $\left|\alpha\right|^2=\tau^2+\lambda$, as well as with $v=e^{-i\rho\cdot x}u$.
With these choices we have
\[\rho\cdot\rho=\lambda,\quad\left|\rho\right|\asymp_\lambda\tau,
\quad\text{and}\quad
\left|e^{-i\rho\cdot x}\right|=e^{-\tau x_n},\]
as well as
\[\left(-\Delta-2i\rho\cdot\nabla\right)v=e^{-i\rho\cdot x}\left(-\Delta-\lambda\right)u.\]

Now Theorem \ref{resolvent-estimate} clearly says that
\[\bigl\|e^{-\tau x_n}\left\langle\cdot\right\rangle^\delta u\bigr\|_{L^2(\mathbb R^n)}
\ll_{\lambda,\delta}\frac1\tau\bigl\|e^{-\tau x_n}\left\langle\cdot\right\rangle^{\delta+1}\left(-\Delta-\lambda\right)u\bigr\|_{L^2(\mathbb R^n)}.\]

\paragraph{Proof of Theorem \ref{superexponential-main-theorem}.}
Let us first observe that by Theorem \ref{exponential-rellich} we must have $\partial^\alpha u\in e^{-\gamma\left\langle\cdot\right\rangle}L^2(\mathbb R^n)$ for each multi-index $\alpha$ with $\left|\alpha\right|\leqslant2$. In particular, we also have
$e^{-\tau x_n}\left\langle\cdot\right\rangle^{\delta+1}\partial^\alpha u\in L^2(\mathbb R^n)$ for all $\tau\in\mathbb R$, any fixed $\delta\in\left]-1,0\right[$ and each multi-index $\alpha$ with $\left|\alpha\right|\leqslant2$.

Our goal will be to prove that $u$ vanishes in $\mathbb R^{n-1}\times\left]-\infty,-2\right[$. Then $u$ also vanishes in $\mathbb R^{n-1}\times\mathbb R_-$ by the unique continuation property of solutions to the free Helmholtz equation.

We want to focus on the behaviour of $u$ in the lower half-space and so we will pick a cut-off function $\chi\in C^{\infty}(\mathbb R^n)$ which only depends on $x_n$, vanishes for $x_n>-1$, and is identically equal to $1$ when $x_n<-2$.

Now Theorem \ref{carleman-estimate} tells us that for large $\tau\in\mathbb R_+$,
\begin{align*}
&e^{2\tau}\bigl\|\left\langle\cdot\right\rangle^\delta u\bigr\|_{L^2(\mathbb R^{n-1}\times\left]-\infty,-2\right[)}\\
&\quad\ll\bigl\|e^{-\tau x_n}\left\langle\cdot\right\rangle^\delta\chi u\bigr\|_{L^2(\mathbb R^n)}\\
&\quad\ll_{\lambda,\delta}\frac1\tau\bigl\|e^{-\tau x_n}\left\langle\cdot\right\rangle^{\delta+1}\left(-\Delta-\lambda\right)(\chi u)\bigr\|_{L^2(\mathbb R^n)}\\
&\quad
=\frac1\tau\left\|e^{-\tau x_n}\left\langle\cdot\right\rangle^{\delta+1}\left(2\,\frac{\partial\chi}{\partial x_n}\cdot\frac{\partial u}{\partial x_n}+\frac{\partial^2\chi}{\partial x_n^2}\,u\right)\right\|_{L^2(\mathbb R^{n-1}\times\left]-2,-1\right[)}\\
&\quad\ll\frac1\tau e^{2\tau}\left\|\left\langle\cdot\right\rangle^{\delta+1}\left(2\,\frac{\partial\chi}{\partial x_n}\cdot\frac{\partial u}{\partial x_n}+\frac{\partial^2\chi}{\partial x_n^2}\,u\right)\right\|_{L^2(\mathbb R^{n-1}\times\left]-2,-1\right[)}
%&\quad\ll_\chi\frac1\tau e^{2\tau}\left\|\left\langle\cdot\right\rangle^{\delta+1}\frac{\partial u}{\partial x_n}\right\|_{L^2(\mathbb R^{n-1}\times\left]-2,-1\right[)}+\frac1\tau e^{2\tau}\bigl\|\left\langle\cdot\right\rangle^{\delta+1}u\bigr\|_{L^2(\mathbb R^{n-1}\times\left]-2,-1\right[)},
\end{align*}
and the result follows by dividing by $e^{2\tau}$ and letting $\tau\longrightarrow\infty$.

\section{Proofs of Theorems \ref{exponential-main-theorem} and \ref{polynomial-main-theorem}}

\subsection{Differentiation under integral signs}

The following lemmas have been adapted from the first chapter of Wong's textbook \cite{Wong}. The reason for them is that we apply analytic continuations for Fourier transforms, and to obtain these extensions with only polynomial decay, we will differentiate under the integral sign in order to check that the Cauchy--Riemann equations hold.
\begin{lemma}
Let $D$ and $\Omega$ be open subsets of $\mathbb R^n$ and $\mathbb R^m$, respectively, and let $f\colon D\times\Omega\longrightarrow\mathbb C$ be measurable. Suppose that
\begin{itemize}\setlength{\itemsep}{0pt}
\item[\rm I.] $f(x,\cdot)\in L^1(\Omega)$ for each $x\in D$;
\item[\rm II.] $f(\cdot,y)\in C^2(D)$ for almost every $y\in\Omega$; and that
\item[\rm III.] for multi-indices $\alpha$ with $\left|\alpha\right|\leqslant2$, we have
\[\sup_{x\in D}\int\limits_\Omega\left|\partial_x^\alpha f(x,y)\right|\mathrm dy<\infty.\]
\end{itemize}
Then, for each $\ell\in\left\{1,2,\ldots,n\right\}$, the integrals
\[\int\limits_\Omega\bigl(\partial_{x_\ell}f\bigr)(x,y)\,\mathrm dy\quad\text{and}\quad
\int\limits_\Omega\left|\bigl(\partial_{x_\ell}f\bigr)(x,y)\right|\mathrm dy\]
are uniformly continuous functions of $x\in D$.
\end{lemma}

\paragraph{Proof.} For each $k\in\left\{1,2,\ldots,n\right\}$ and small $h\in\mathbb R$, we may apply the mean value theorem to estimate
\begin{align*}
&\int\limits_\Omega\left|
\bigl(\partial_{x_\ell}f\bigr)(x_1,\ldots,x_k+h,\ldots,x_n,y)
-\bigl(\partial_{x_\ell}f\bigr)(x_1,\ldots,x_k,\ldots,x_n,y)\right|\mathrm dy\\
&=\left|h\right|\int\limits_\Omega\left|\bigl(\partial_{x_k}\partial_{x_\ell}f\bigr)(x_1,\ldots,\xi(x,y;h),\ldots,x_n,y)\right|\mathrm dy\\
&\leqslant\left|h\right|\sup_{x\in D}\int\limits_\Omega\left|\bigl(\partial_{x_k}\partial_{x_\ell}f\bigr)(x,y)\right|\mathrm dy.
\end{align*}
Here $\xi(\ldots)$ has the obvious meaning.

\begin{theorem}\label{differentiation-under-integral-sign}
Let $D$ and $\Omega$ be open subsets of $\mathbb R^n$ and $\mathbb R^m$, respectively, and let $f\colon D\times\Omega\longrightarrow\mathbb C$ be measurable. Suppose that
\begin{itemize}\setlength{\itemsep}{0pt}
\item[\rm I.] $f(x,\cdot)\in L^1(\Omega)$ for each $x\in D$;
\item[\rm II.] $f(\cdot,y)\in C^2(D)$ for almost every $y\in\Omega$; and that
\item[\rm III.] for all multi-indices $\alpha$ with $\left|\alpha\right|\leqslant2$, we have
\[\sup_{x\in D}\int\limits_\Omega\left|\bigl(\partial_x^\alpha f\bigr)(x,y)\right|\mathrm dy<\infty.\]
\end{itemize}
Then
\[\int\limits_\Omega f(\cdot,y)\,\mathrm dy\in C^1(D)\]
and
\[\partial_{x_\ell}\int\limits_\Omega f(x,y)\,\mathrm dy
=\int\limits_\Omega\bigl(\partial_{x_\ell}f\bigr)(x,y)\,\mathrm dy\]
in $D$ for each $\ell\in\left\{1,2,\ldots,n\right\}$.
\end{theorem}

\paragraph{Proof.}
By the fundamental theorem of analysis, the previous lemma and Fubini's theorem,
\begin{align*}
&\partial_{x_\ell}\int\limits_\Omega f(x,y)\,\mathrm dy\\
&=\lim_{h\longrightarrow0}\frac1h\int\limits_\Omega
\bigl(f(x_1,\ldots,x_\ell+h,\ldots,x_n,y)-f(x_1,\ldots,x_\ell,\ldots,x_n,y)\bigr)\,\mathrm dy\\
&=\lim_{h\longrightarrow0}\frac1h\int\limits_\Omega
\int\limits_{x_\ell}^{x_\ell+h}\bigl(\partial_{x_\ell}f\bigr)(x_1,\ldots,x_{\ell-1},s,x_{\ell+1},\ldots,x_n,y)\,\mathrm ds\,\mathrm dy\\
&=\lim_{h\longrightarrow0}\frac1h
\int\limits_{x_\ell}^{x_\ell+h}\int\limits_\Omega
\bigl(\partial_{x_\ell}f\bigr)(x_1,\ldots,x_{\ell-1},s,x_{\ell+1},\ldots,x_n,y)\,\mathrm dy\,\mathrm ds\\
&=\int\limits_\Omega\bigl(\partial_{x_\ell}f\bigr)(x,y)\,\mathrm dy.
\end{align*}

\subsection{Some Paley--Wiener theorems% and a Rellich type theorem for exponential decay
}

The following classical result is e.g.\ Theorem 19.2 in Rudin's textbook \cite{Rudin}.
\begin{theorem}\label{half-line-paley-wiener}
A function $f\in L^2(\mathbb R)$ is supported in $\left[0,\infty\right[$ \textbf{if and only if} its Fourier transform $\widehat f$ extends to an analytic function in the lower half-plane $\left\{z\in\mathbb C\middle|\Im z<0\right\}$ and this continuation satisfies
\[\sup_{\eta\in\mathbb R_-}\bigl\|\widehat f(\cdot+i\eta)\bigr\|_{L^2(\mathbb R)}<\infty.\]
\end{theorem}
\noindent The following Paley--Wiener theorem is e.g.\ Theorem XI.13 in \cite[p.~18]{Reed--SimonII}.
\begin{theorem}\label{exponential-paley-wiener}
Let $\gamma_0\in\mathbb R_+$, and let $f\in L^2(\mathbb R^n)$. Then $e^{\gamma\left\langle\cdot\right\rangle}\,f\in L^2(\mathbb R^n)$ for each $\gamma\in\left]0,\gamma_0\right[$ \textbf{if and only if}
the Fourier transform $\widehat f\in L^2(\mathbb R^n)$ extends analytically to the set $\left\{\zeta\in\mathbb C^n\middle|\left|\Im\zeta\right|<\gamma_0\right\}$ so that, for each $\eta\in\mathbb R^n$ with $\left|\eta\right|<\gamma_0$, we have $\widehat f(\cdot+i\eta)\in L^2(\mathbb R^n)$, and that
for each $\gamma\in\left]0,\gamma_0\right[$,
\[\sup_{\substack{\eta\in\mathbb R^n\\\left|\eta\right|\leqslant\gamma}}
\bigl\|\widehat f(\cdot+i\eta)\bigr\|_{L^2(\mathbb R^n)}<\infty.\]
\end{theorem}

%The following is essentially Lemma 2.5 in \cite[p.~1780]{Paivarinta--Salo--Uhlmann}, which may be reformulated in the following way in view of Theorem 14.3.6 in \cite{Hormander-II}.
%\begin{theorem}\label{exponential-rellich}
%Let $\gamma_0\in\mathbb R_+$, and let $f$ be a function such that $e^{\gamma\left\langle\cdot\right\rangle}\,f\in L^2(\mathbb R^n)$ for each $\gamma\in\left]0,\gamma_0\right[$. If $u\in\mathring B^\ast_2$ solves the equation
%\[\left(-\Delta-\lambda\right)u=f,\]
%then $e^{\gamma\left\langle\cdot\right\rangle}\,u\in H^2(\mathbb R^n)$ for each $\gamma\in\left]0,\gamma_0\right[$.
%\end{theorem}
%\noindent For us, the relevant value will be $\gamma=\frac{\gamma_0}2$ when proving the discreteness of non-scattering energies for exponentially decaying potentials.

\subsection{Division by the symbol on the Fourier side}

For the rest of this section we will simplify our notation by writing $p$ for the symbol polynomial $4\pi^2\!\left(z_1^2+z_2^2+\ldots+z_n^2\right)$. We shall also consider the level-set manifolds
\[M_\lambda^\mathbb R=\left\{\xi\in\mathbb R^n\,\middle|\,p\!\left(\xi\right)=\lambda\right\},\]
and
\[M_\lambda^\mathbb C=\left\{\zeta\in\mathbb C^n\,\middle|\,p\!\left(\zeta\right)=\lambda\right\},\]
where $\lambda\in\mathbb R_+$, as usual.
%We let $F'$ denote the Fourier transform with respect to the first $n-1$ variables, and $F_n$ with respect to the last variable.

\begin{lemma}\label{lemma-from-PSU}
Let $D\subseteq\mathbb C^n$ be an open set such that $M_\lambda^{\mathbb R}\subseteq D$ and $M_\lambda^{\mathbb C}\cap D$ is connected. If $f\colon D\longrightarrow\mathbb C$ is analytic and vanishes on the real sphere $M_\lambda^{\mathbb R}$,
then $f$ also vanishes in the intersection $M_\lambda^{\mathbb C}\cap D$ and the expression $f/(p-\lambda)$ gives rise to an analytic function in $D$.
\end{lemma}
\noindent
The proof is modelled after a portion of the proof of Lemma 2.5 in \cite{Paivarinta--Salo--Uhlmann}.

\paragraph{Proof.}
Let $\xi\in M_\lambda^{\mathbb R}$ be arbitrary. At least one coordinate of $\xi$ must be non-zero, say $\xi_n\neq0$. Thus, we have $\nabla p=8\pi^2\xi\neq0$ at $\xi$, and by the inverse function theorem, the mapping
\[\varphi=\zeta\longmapsto\left\langle\zeta',p\!\left(\zeta\right)-\lambda\right\rangle\colon W\longrightarrow\mathbb C^n\]
is a biholomorphic diffeomorphism between an open connected neighbourhood $W$ of $\xi$ in $D$ and $\varphi\!\left[W\right]$. This particular mapping ``straightens out'' a portion of the sphere around $\xi$:
\[\varphi\bigl[M_\lambda^{\mathbb C}\cap W\bigr]
=\bigl\{\zeta\in\varphi\!\left[W\right]\bigm|\zeta_n=0\bigr\}.\]

The function $f\circ\varphi^{-1}\colon\varphi\bigl[W\bigr]\longrightarrow\mathbb C$ vanishes on $\varphi\bigl[W\bigr]\cap\left(\mathbb R^{n-1}\times\left\{0\right\}\right)$. Therefore it also vanishes on $\varphi\bigl[W\bigr]\cap\left(\mathbb C^{n-1}\times\left\{0\right\}\right)$, which in turn implies that
\[f\big|_W=f\circ\varphi^{-1}\circ\varphi\] vanishes in $M_\lambda^{\mathbb C}\cap W$.

Now $f$ vanishes in a neighbourhood of $M_\lambda^{\mathbb R}$ in $M_\lambda^{\mathbb C}\cap D$ and must vanish everywhere on $M_\lambda^{\mathbb C}\cap D$ as an analytic function on a connected analytic manifold.

The second part of the lemma is proved by a similar reasoning. We pick again an arbitrary point $\xi$, but this time from the complex sphere $M_\lambda^{\mathbb C}\cap D$. Again we will have a biholomorphic diffeomorphism $\varphi$ from some open connected neighbourhood $W$ of $\xi$ in $D$ into $\varphi[W]$, given by the same expression as before. Since $f\circ\varphi^{-1}$ vanishes on $\varphi[W]\cap\left(\mathbb C^{n-1}\times\left\{0\right\}\right)$, the expression $\left(f\circ\varphi^{-1}\right)\!(\zeta)/\zeta_n$ defines a function analytic in $\varphi[W]$. Finally, we conclude that
\[\frac f{p-\lambda}=\left(\zeta\longmapsto\frac{\left(f\circ\varphi^{-1}\right)\left(\zeta\right)}{\zeta_n}\right)\circ\varphi\]
is analytic in $W$, and since $\xi$ was arbitrary, and since $\frac f{p-\lambda}$ is definitely analytic outside of $M_\lambda^{\mathbb C}$, we are done.

\subsection{Proving Theorems \ref{exponential-main-theorem} and \ref{polynomial-main-theorem}}

The following proof works verbatim for both Theorems \ref{exponential-main-theorem} and \ref{polynomial-main-theorem}. The only difference is in the reasons for the analytic continuations of $\widehat f$.
The proof is modelled after the proof of Lemma 2.5 in \cite{Paivarinta--Salo--Uhlmann} and the proof of Theorem 8.3 in \cite{Hitrik--Krupchyk--Ola--Paivarinta}.

Taking Fourier transforms of both sides of the Helmholtz equation gives
\[\left(p-\lambda\right)\widehat u=\widehat f,\]
which holds in $\mathbb R^n$.
A basic result in scattering theory, Theorem 14.3.6 from \cite{Hormander-II}, says that $\widehat f\,\big|_{M_\lambda^{\mathbb R}}\equiv0$.

The assumptions on $f$ guarantee that $\widehat f$ extends to an analytic function in
\[D=\bigl\{\zeta\in\mathbb C^n\bigm|\left|\Im\zeta\right|<\gamma_0\bigr\},\]
where the constant $\gamma_0\in\mathbb R_+$ is the same as in the assumptions of the theorem. For Theorems \ref{exponential-main-theorem}, this follows immediately from Theorem \ref{exponential-paley-wiener}.
For Theorem \ref{polynomial-main-theorem} we use Theorem \ref{differentiation-under-integral-sign} which guarantees that the Fourier transform $\widehat f$ can be differentiated under the integral sign and so, from the Cauchy--Riemann equations for the integrand, we get the Cauchy--Riemann equations for $\widehat f$.

In any case, combining the above facts with Lemma \ref{lemma-from-PSU} leads to the conclusion that the expression $\widehat f/(p-\lambda)$ gives rise to an analytic function in~$D$. In particular, $\widehat u$ has an analytic extension to~$D$.

Next, let us fix a point $\xi'\in B^{n-1}(0,\sqrt\lambda/2\pi)\subseteq\mathbb R^{n-1}$. For clarity, we write $q\!\left(\cdot\right)$ for $p\!\left(\xi',\cdot\right)$. Then $q-\lambda$ is an entire function of one complex variable, and its only zeroes are simple ones at the points
\[\pm\mu=\pm\frac1{2\pi}\sqrt{\lambda-4\pi^2\left|\xi'\right|^2}.\]

Since $f$ vanishes in $\mathbb R^{n-1}\times\mathbb R_-$, the Fourier transform $F'f(\xi',\cdot)$ vanishes in $\mathbb R_-$, so that by Theorem \ref{half-line-paley-wiener} $\widehat f$ has an analytic extension in the last variable to $\mathbb R\times i\left]-\infty,\gamma_0\right[$, and
\[\int\limits_{-\infty}^\infty\bigl|\widehat f\!\left(\xi',\xi_n-i\eta\right)\bigr|^2\,\mathrm d\xi_n\ll\int\limits_{-\infty}^\infty\bigl|\widehat f\!\left(\xi',\xi_n\right)\bigr|^2\,\mathrm d\xi_n<\infty\]
for all $\eta\in\mathbb R_+$. Of course, $\widehat u(\xi',\cdot)$ has an analytic extension to $\mathbb R\times i\left]-\infty,\gamma_0\right[$ as well.

Since $\left|q\!\left(z\right)-\lambda\right|$ is bounded from below, when $z\in\mathbb C$ and $\Im z<-1$, we have
\begin{align*}
\int\limits_{-\infty}^\infty\bigl|\widehat u\!\left(\xi',\xi_n-i\eta\right)\bigr|^2\,\mathrm d\xi_n
&=\int\limits_{-\infty}^\infty\left|\frac{\widehat f\!\left(\xi',\xi_n-i\eta\right)}{q\!\left(\xi_n-i\eta\right)-\lambda}\right|^2\,\mathrm d\xi_n\ll\int\limits_{-\infty}^\infty\bigl|\widehat f\!\left(\xi',\xi_n\right)\bigr|^2\,\mathrm d\xi_n,
\end{align*}
whenever $\eta\in\left[1,\infty\right[$. In the same vein, the expression
\[\left(\,\int\limits_{-\infty}^{-2\mu}+\int\limits_{2\mu}^\infty\,\right)
\bigl|\widehat u\!\left(\xi',\xi_n-i\eta\right)\bigr|^2\,\mathrm d\xi_n\]
is also bounded by a constant independent of $\eta$,
whenever $\eta\in\left[0,1\right[$.

Also, since the function $\widehat f(\xi',\cdot)/(q-\lambda)$ is analytic in a neighbourhood of the rectangle $\left[-2\mu,2\mu\right]\times i\left[-1,0\right]\subseteq\mathbb C$, it is bounded as well, and so
\[\int\limits_{-2\mu}^{2\mu}\bigl|\widehat u\!\left(\xi',\xi_n-i\eta\right)\bigr|^2\,\mathrm d\xi_n\]
is bounded from above by something constant and independent from $\eta$, even for $\eta\in\left[0,1\right[$.

Now we are able to conclude from Theorem \ref{half-line-paley-wiener} that
\[F'u\!\left(\xi',x_n\right)=F_n^{-1}\widehat u\!\left(\xi',x_n\right)=0\]
for all $\xi'\in B^{n-1}\!\left(0,\frac{\sqrt\lambda}{2\pi}\right)$ and $x_n\in\mathbb R_-$. Since $F'u$ is analytic with respect to the first $n-1$ variables, $F'u\!\left(\xi',x_n\right)=0$ for all $\xi'\in\mathbb R^{n-1}$ and all $x_n\in\mathbb R_-$. Finally, taking $F'^{-1}$ gives the desired conclusion that $u$ vanishes in $\mathbb R^{n-1}\times\mathbb R_-$.

\section{Proof of Theorem~\ref{hard-corollary}}\label{corollary-section}

Let $\lambda$ be a non-scattering energy for $V$. Then
\[\left\{\!\!\begin{array}{l}
\left(-\Delta+V-\lambda\right)v=0,\\
\left(-\Delta-\lambda\right)w=0
\end{array}\right.\]
in $\mathbb R^n$ for some $v,w\in B_2^\ast\setminus0$ with $v-w\in\mathring B_2^\ast$. Now consider $u=v-w$, which solves
\[\left(-\Delta-\lambda\right)u=-Vv.\]
By Theorem~\ref{exponential-main-theorem} and the assumptions of Theorem~\ref{hard-corollary}, the function $u$ vanishes in $\mathbb R^n\setminus\Omega$. Furthermore,
Theorem~\ref{exponential-rellich} says that $u$ belongs to the space $H^2_0(\Omega;e^{\gamma\left\langle\cdot\right\rangle})$ with $\gamma=\gamma_0/2$, which is the closure of test functions $u\in C_{\mathrm c}^\infty(\Omega)$ with respect to the weighted Sobolev norm
\[\bigl\|e^{\gamma\left\langle\cdot\right\rangle}u\bigr\|
+\bigl\|e^{\gamma\left\langle\cdot\right\rangle}\nabla u\bigr\|
+\bigl\|e^{\gamma\left\langle\cdot\right\rangle}\nabla\otimes\nabla u\bigr\|,\]
where $\left\|\cdot\right\|$ denotes the usual $L^2$-norm.

We shall use $H_0^2(\Omega;e^{\gamma\left\langle\cdot\right\rangle})$ as a quadratic form domain. As the ambient Hilbert space we shall use the space $L^2(\Omega;e^{\gamma\left\langle\cdot\right\rangle})$, defined
in the obvious way by
the weighted norm $\left\|e^{\gamma\left\langle\cdot\right\rangle}u\right\|$.

Let us consider the composition of mappings
\[H_0^2(\Omega;e^{\gamma\left\langle\cdot\right\rangle})
\longrightarrow H_0^2(\Omega)
\longrightarrow L^2(\Omega)
\longrightarrow L^2(\Omega;e^{\gamma\left\langle\cdot\right\rangle}),\]
where the middle mapping is the compact embedding, and the first and the last mappings are just multiplication and division by $e^{\gamma\left\langle\cdot\right\rangle}$, respectively. We easily see that the first and last mappings are bounded, and so $H_0^2(\Omega;e^{\gamma\left\langle\cdot\right\rangle})$ embeds compactly into $L^2(\Omega;e^{\gamma\left\langle\cdot\right\rangle})$.

We have now reduced the situation to a single fourth-order equation:
\begin{lemma}
Under the assumptions of Theorem~\ref{hard-corollary}, if $\lambda\in\mathbb R_+$ is a non-scattering energy, then there exists a function $u\in H_0^2(\Omega;e^{\gamma\left\langle\cdot\right\rangle})\setminus0$ solving the fourth-order equation
\begin{equation}\label{fourth-order-equation}
\left(-\Delta+V-\lambda\right)\frac1V\left(-\Delta-\lambda\right)u=0
\end{equation}
in $\Omega$ in the sense of distributions. Furthermore, this transition respects multiplicities.
\end{lemma}
\noindent The discreteness of non-scattering energies will therefore follow from the following proposition.
\begin{theorem}\label{discrete-fourth-order}
The set of real numbers $\lambda$ for which the equation (\ref{fourth-order-equation}) has a non-trivial $H_0^2(\Omega;e^{\gamma\left\langle\cdot\right\rangle})$-solution is a discrete subset of $\left[0,\infty\right[$. For each such $\lambda$ the space of solutions is finite dimensional. \end{theorem}

The operator on the left-hand side of (\ref{fourth-order-equation}) can be treated nicely via quadratic forms, and for this purpose we define for each $\lambda\in\mathbb C$ the quadratic
form
\[Q_\lambda=u\longmapsto\left\langle(-\Delta+V-\overline\lambda)\,u\,\middle|\,\frac1V\,
(-\Delta-\lambda)\,u\right\rangle\colon H_0^2(\Omega;e^{\gamma\left\langle\cdot\right\rangle})\longrightarrow\mathbb C,\]
where the $L^2$-inner product is linear in the second argument.
The family $\langle
Q_\lambda\rangle_{\lambda\in\mathbb C}$ has the pleasant
properties enumerated in the theorem below. These properties are analogous to a part of Theorem 4 of \cite{Vesalainen}.
\begin{theorem}\label{list-of-good-properties}
\begin{enumerate}
\item
The quadratic forms $Q_\lambda$ form an entire self-adjoint analytic family of forms of type (a) with compact resolvent, and therefore gives rise to a family of operators $T_\lambda$, which is an entire self-adjoint analytic family of operators of type (B) with compact resolvent.

\item
Furthermore, there exists a sequence $\left\langle\mu_\nu(\cdot)\right\rangle_{\nu=1}^\infty$ of real-analytic functions $\mu_\nu(\cdot)\colon\mathbb R\longrightarrow\mathbb R$ such that, for real $\lambda$, the spectrum of $T_\lambda$, which consists of a discrete set of real eigenvalues of finite multiplicity, consists of $\mu_1(\lambda)$, $\mu_2(\lambda)$, \dots, including multiplicity.

\item
In addition, for any given $T\in\mathbb R_+$, there exists constant $c\in\mathbb R_+$ such that
\[\bigl|\mu_\nu(\lambda)-\mu_\nu(0)\bigr|\ll_Te^{c|\lambda|}-1\]
for all $\lambda\in\left[-T,T\right]$ and each $\nu\in\mathbb Z_+$.

\item
The pairs $\left\langle\lambda,u\right\rangle\in\mathbb R\times H_0^2(\Omega;e^{\gamma\left\langle\cdot\right\rangle})$ for which (\ref{fourth-order-equation}) holds,
are in bijective correspondence with the pairs $\left\langle\nu,\lambda\right\rangle\in\mathbb Z_+\times\mathbb R$ for which $\mu_\nu(\lambda)=0$.
\end{enumerate}
\end{theorem}

Theorem \ref{discrete-fourth-order} follows easily from these properties of $Q_\lambda$: It is obvious
that zero is not an eigenvalue of $T_\lambda$ for any negative real $\lambda$,
as $Q_\lambda(u)>0$ for all non-zero functions $u\in\mathrm{Dom}\,Q_\lambda$. Hence none
of the functions $\mu_\nu(\cdot)$ can vanish identically, so that the set of
zeroes of each of them is discrete. Why the union of the zero sets can not have
an accumulation point follows immediately
from the third statement above,
which says that the functions $\mu_\nu(\cdot)$ change their values uniformly
locally exponentially. That is, when the value of $\lambda$ changes by a finite amount, only finitely many $\mu_\nu(\cdot)$ will have enough time to drop to zero, and the discreteness has been obtained.

The proof of Theorem~\ref{list-of-good-properties} depends heavily on the basic theory of quadratic forms and analytic perturbation theory. For an excellent reference on these topics, we recommend the book by Kato \cite{Kato}, in particular its Chapters VI and VII.

As the arguments in \cite{Vesalainen}, the proof of Theorem 4 there to be precise, work verbatim in our case, except for the required weighted inequality, which is given below,  we simply refer the reader to \cite{Vesalainen}.
The following weighted inequality replaces Lemma 2 of \cite{Vesalainen}.
\begin{lemma}
Let $\gamma\in\mathbb R_+$ with $\gamma\ll_n1$, and let us be given a compact subset $K\subseteq\mathbb C$. Then we have, for all $\lambda\in K$ and all $u\in C_{\mathrm c}^\infty(\mathbb R^n)$, the weighted inequality
\[\bigl\|e^{\gamma\left\langle\cdot\right\rangle}\,u\bigr\|
+\bigl\|e^{\gamma\left\langle\cdot\right\rangle}\,\nabla u\bigr\|
+\bigl\|e^{\gamma\left\langle\cdot\right\rangle}\,\nabla\otimes\nabla u\bigr\|
\ll_{n,K}
\bigl\|e^{\gamma\left\langle\cdot\right\rangle}\left(-\Delta-\lambda\right)u\bigr\|
+\bigl\|e^{\gamma\left\langle\cdot\right\rangle}\,u\bigr\|.\]
\end{lemma}

\paragraph{Proof.}
The elementary inequalities
\[\left\langle\cdot\right\rangle^4\ll\bigl|4\pi^2\left|\cdot\right|^2+1\bigr|^2
\ll\bigl|4\pi^2\left|\cdot\right|^2-\lambda\bigr|^2+\left|\lambda+1\right|^2\]
imply that
\[\bigl\|u\bigr\|+\bigl\|\nabla u\bigr\|+\bigl\|\nabla\otimes\nabla u\bigr\|
\ll_n\bigl\|\left(-\Delta-\lambda\right)u\bigr\|
+\left\langle\lambda\right\rangle\bigl\|u\bigr\|.\]
Now we can introduce the exponential weights into this applying Leibniz's rule:
\begin{align*}
&\bigl\|e^{\gamma\left\langle\cdot\right\rangle}\,u\bigr\|
+\bigl\|e^{\gamma\left\langle\cdot\right\rangle}\,\nabla u\bigr\|
+\bigl\|e^{\gamma\left\langle\cdot\right\rangle}\,\nabla\otimes\nabla u\bigr\|\\
&\qquad\ll_n
\bigl\|e^{\gamma\left\langle\cdot\right\rangle}\,u\bigr\|
+\bigl\|\nabla\bigl(e^{\gamma\left\langle\cdot\right\rangle}\,u\bigr)\bigr\|
+\bigl\|\nabla\otimes\nabla\bigl(e^{\gamma\left\langle\cdot\right\rangle}\,u\bigr)\bigr\|\\
&\qquad\qquad+\bigl\|\bigl(\nabla e^{\gamma\left\langle\cdot\right\rangle}\bigr)\,u\bigr\|
+\bigl\|\bigl(\nabla e^{\gamma\left\langle\cdot\right\rangle}\bigr)\otimes\nabla u\bigr\|
+\bigl\|\bigl(\nabla\otimes\nabla e^{\gamma\left\langle\cdot\right\rangle}\bigr)\,u\bigr\|\\
&\qquad\ll_n
\bigl\|\left(-\Delta-\lambda\right)\bigl(e^{\gamma\left\langle\cdot\right\rangle}\,u\bigr)\bigr\|
+\left\langle\lambda\right\rangle\bigl\|e^{\gamma\left\langle\cdot\right\rangle}\,u\bigr\|\\
&\qquad\qquad
+\gamma\bigl\|e^{\gamma\left\langle\cdot\right\rangle}\,u\bigr\|
+\gamma\bigl\|e^{\gamma\left\langle\cdot\right\rangle}\,\nabla u\bigr\|
+\left(\gamma+\gamma^2\right)\bigl\|e^{\gamma\left\langle\cdot\right\rangle}\,u\bigr\|\\
&\qquad\ll_n
\bigl\|e^{\gamma\left\langle\cdot\right\rangle}\left(-\Delta-\lambda\right)u\bigr\|
+\bigl\|\bigl(\nabla e^{\gamma\left\langle\cdot\right\rangle}\bigr)\cdot\nabla u\bigr\|
+\bigl\|\bigl(\Delta e^{\gamma\left\langle\cdot\right\rangle}\bigr)\,u\bigr\|\\
&\qquad\qquad
+\left(\left\langle\lambda\right\rangle+\gamma+\gamma^2\right)\bigl\|e^{\gamma\left\langle\cdot\right\rangle}\,u\bigr\|
+\gamma\bigl\|e^{\gamma\left\langle\cdot\right\rangle}\,\nabla u\bigr\|\\
&\qquad\ll_n
\bigl\|e^{\gamma\left\langle\cdot\right\rangle}\left(-\Delta-\lambda\right)u\bigr\|
+\left(\left\langle\lambda\right\rangle+\gamma+\gamma^2\right)\bigl\|e^{\gamma\left\langle\cdot\right\rangle}\,u\bigr\|
+\gamma\bigl\|e^{\gamma\left\langle\cdot\right\rangle}\,\nabla u\bigr\|.
\end{align*}
Finally, the last term may be absorbed to the original left-hand side provided that $\gamma\ll_n1$.

\section{Proof of Theorem \ref{discrete-main-theorem}}

We begin with the following analogue of the Paley--Wiener theorem on exponential decay of the Fourier transform for functions defined in $\mathbb Z^n$. It could be compared to, say, Theorem IX.13 of \cite{Reed--SimonII}.
\begin{theorem}\label{discrete-exponential-paley-wiener}
Let $f\in\ell^2(\mathbb Z^n)$ and let $\gamma_0\in\mathbb R_+$. Then $e^{\gamma\left\langle\cdot\right\rangle}f\in\ell^2(\mathbb Z^n)$ for all $\gamma\in\left]0,\gamma_0\right[$ \textbf{if and only if}
the function $\check f\in L^2(\mathbb T^n)$ extends to an analytic function in \[\bigl\{z\in\mathbb T^n_{\mathbb C}\bigm|\left|\Im z\right|<\gamma_0/(2\pi)\bigr\}.\]
\end{theorem}

\paragraph{Proof.}
First, let $\gamma_0\in\mathbb R$ and $f\in\ell^2(\mathbb Z^n)$ be such that $e^{\gamma\left\langle\cdot\right\rangle}f\in\ell^2(\mathbb Z^n)$ for all $\gamma\in\left]0,\gamma_0\right[$. Fix some $\gamma\in\left]0,\gamma_0\right[$. Then certainly
$e^{\varepsilon\left\langle\cdot\right\rangle}e^{2\pi\eta\cdot\xi}f(x)\in\ell^2(\mathbb Z^n)$ for any $\eta\in\mathbb R^n$ with $\left|\eta\right|\leqslant\gamma/2\pi$, for small $\varepsilon\in\mathbb R_+$, and the series
\[\check f(z)=\sum_{\xi\in\mathbb Z^n}f(\xi)\,e(\xi\cdot z)\]
clearly converges absolutely and uniformly for all $z\in\mathbb T_{\mathbb C}^n$ with $\left|\Im z\right|\leqslant\gamma/2\pi$, and this limit must be analytic in $z$ since each of the terms is.

Next, assume that $f\in\ell^2(\mathbb Z^n)$ and $\gamma_0\in\mathbb R_+$ are such that $\check f\in L^2(\mathbb R^n)$ extends to an analytic function in $\left\{z\in\mathbb T_{\mathbb C}^n\middle|\left|\Im z\right|<\gamma_0/2\pi\right\}$. Fix some $\gamma\in\left]0,\gamma_0\right[$. Now the restriction of $|\check f|$ to the compact set $\left\{z\in\mathbb T_{\mathbb C}^n\middle|\left|\Im z\right|\leqslant\gamma/2\pi\right\}$ must be uniformly bounded by some constant $C_\gamma\in\mathbb R_+$ only depending on $\gamma$.

Then, for arbitrary $y\in\mathbb R^n$ with $\left|y\right|=\gamma/2\pi$, we may estimate, using Cauchy's integral theorem, that for any given $\xi\in\mathbb Z^n$,
\begin{align*}
f(\xi)
&=\int\limits_{\mathbb T^n}\check f(x)\,e(-x\cdot\xi)\,\mathrm dx\\
&=\int\limits_{\mathbb T^n}\check f(x-iy)\,e\bigl(-(x-iy)\cdot\xi)\bigr)\,\mathrm dx\\
&\ll e^{2\pi y\cdot\xi}\int\limits_{\mathbb T^n}\bigl|\check f(x-iy)\bigr|\,\mathrm dx\ll e^{2\pi y\cdot\xi}\,C_\gamma.
\end{align*}
Thus, we have
\[f(\xi)\ll_\gamma\inf_{\substack{y\in\mathbb R^n,\\\left|y\right|=\gamma/2\pi}}e^{2\pi y\cdot\xi}=e^{-\gamma\left|\xi\right|}.\]

\bigbreak
The following is, more or less, a discrete analogue of
Theorem \ref{exponential-rellich}.
%Lemma 2.5 of \cite{Paivarinta--Salo--Uhlmann}.
\begin{theorem}\label{discrete-exponential-rellich}
Let $f\in\ell^2(\mathbb Z^n)$ be such that $e^{\gamma\left\langle\cdot\right\rangle}f\in\ell^2(\mathbb Z^n)$ for all $\gamma\in\mathbb R_+$. Also, let $u\colon\mathbb Z^n\longrightarrow\mathbb C$ be such that
\[\frac1R\sum_{\left|\xi\right|\leqslant R}\bigl|u(\xi)\bigr|^2\longrightarrow0\]
as $R\longrightarrow\infty$. Finally, let $\lambda\in\left]0,n\right[$, and assume that
\[\left(-\Delta_{\mathrm{disc}}-\lambda\right)u=f\]
in $\mathbb Z^n$. Then also $e^{\gamma\left\langle\cdot\right\rangle}u\in\ell^2(\mathbb Z^n)$ for all $\gamma\in\mathbb R_+$.
\end{theorem}

This follows easily: By Theorem \ref{discrete-exponential-paley-wiener} the Fourier series $\check f$ extends to an entire function in $\mathbb T^n_{\mathbb C}$. Furthermore,
\[\bigl(h(x)-\lambda\bigr)\,\check u=\check f\]
for $x\in\mathbb T^n$, where
\[h(x)=\sum_{j=1}^n\sin^2\frac{x_j}2.\] 
We point out that even though $u$ doesn't strictly speaking belong to $\ell^2$, it is certainly at most polynomially growing, allowing us to consider $\check u$ as a distribution; for more on this point of view, see e.g. Chapter 3 in the book \cite{Ruzhansky--Turunen}.

Now the arguments of Section 4 of \cite{Isozaki--Morioka} show that $\check u$ extends to an entire function in $\mathbb T^n_{\mathbb C}$. Namely, the arguments in Section 4.1 do not involve $\check f$ at all, and in Section 4.2, the proof of Lemma 4.3 only depends on the smoothness of $\check f$, and after Lemma 4.4, when the analytic continuation of $\check u$ is obtained, only the analytic continuation of $\check f$ is required.
Finally, the conclusion follows from Theorem \ref{discrete-exponential-paley-wiener}.

\paragraph{Proof of Theorem \ref{discrete-main-theorem}.}

We shall prove that $u(0)=0$. Given a point $\xi_0\in C$, the same argument applied to $u(\cdot+\xi_0)$ and $f(\cdot+\xi_0)$ shows that $u(\xi_0)=0$.

The idea is to apply the definition of $\Delta_{\mathrm{disc}}$ and the discrete Helmholtz equation in the form:
\begin{equation}\label{trick}
u(\xi)=\left(2n-4\lambda\right)\,u(\xi+e_n)-u(\xi+2e_n)-\sum_{j=1}^{n-1}\bigl(u(\xi+e_n+e_j)+u(\xi+e_n-e_j)\bigr).
\end{equation}
This holds for all $\xi\in C$.
Applying this once to $u(0)$ gives $\leqslant2n$ terms of the form $u(\xi)$ with $\xi\in C$ and $1\leqslant\xi_n\leqslant2$, with constant coefficients, each of which has absolute value $\leqslant2n$.

Applying \eqref{trick} again to all the terms of the previous step gives rise to $\leqslant4n^2$ terms of the form $u(\xi)$ with $\xi\in C$ and $2\leqslant\xi_n\leqslant4$, with coefficients of size $\leqslant4n^2$.

Continuing in this manner, after $N\in\mathbb Z_+$ steps $u(0)$ has been represented as the sum of $\leqslant(2n)^N$ terms of the form $u(\xi)$ with $\xi\in C$ and $N\leqslant\xi_n\leqslant2N$, with coefficients of size $\leqslant(2n)^N$. Thus, by the triangle inequality,
\[\left|u(0)\right|\leqslant(2n)^N\,(2n)^N\max_{\substack{\xi\in C,\\N\leqslant\xi_n\leqslant2N}}\left|u(\xi)\right|.\]

Theorem \ref{discrete-exponential-rellich} tells us that $u(\xi)\ll_\gamma e^{-\gamma\left\langle\xi\right\rangle}$ for all $\xi\in\mathbb Z^n$ and any fixed $\gamma\in\mathbb R_+$. In particular,
\[u(0)\ll_\gamma (4n^2)^N\,e^{-\gamma N},\]
and choosing $\gamma>\log4n^2$ and letting $N\longrightarrow\infty$ gives the desired result.

\section*{Acknowledgements}

This research was funded by Finland's Ministry of Education through the Doctoral School of Inverse Problems, and by the Finnish Centre of Excellence in Inverse Problems Research.
The author wishes to thank Prof.\ M. Salo for suggesting that Rellich's lemma might admit generalizations for unbounded domains and for pointing out Corollary \ref{carleman-estimate}, as well as for a number of helpful and encouraging discussions.

\end{document}